\documentstyle{amsppt}

\def\ms{\medskip}

\def\varep{\varepsilon}

\def\ms{\medskip}

\def\cR{\Cal R}

\def\varep{\varepsilon}

\def\bbox{\vrule height.6em width.6em depth0em} 


\overfullrule=0pt

 at 10truept

\def\r{\Bbb R}
\def\re{\Bbb R}

\def\t{\Bbb T}
\def\z{\Bbb Z}

\topmatter
\title The Rotation Set and Periodic Points for Torus Homeomorphisms \endtitle
\author by John Franks \endauthor
\date Department of Mathematics, Northwestern University\enddate

\abstract
We consider the rotation set $\rho(F)$ for a lift $F$ of an area
preserving homeomorphism $f: \t^2\to \t^2$, which is homotopic to the
identity. The relationship between this set and the existence of
periodic points for $f$ is least well understood in the case when this
set is a line segment.  We show that in this case if a vector $v$ lies in 
$\rho(F)$ and has both co-ordinates rational, then there is a
periodic point $x\in \t^2$ with the property that
$$\frac{F^q(x_0)-x_0}q = v$$ where $x_0\in \re^2$ is any lift of $x$ and
$q$ is the least period of $x$.
\endabstract
\endtopmatter

\document

One of the simplest and most important invariants for the study of
homeomorphisms of the circle is the so called rotation number introduced
by Poincar\'e.  It is defined in terms of the average behavior of
the iterates of points, but turns out to be independent of the point
used to define it.  Its great utility lies in a close relationship
between the rationality of the rotation number and the existence of
critical points.  

Considerable effort has been devoted to attempts to generalize the
results concerning the rotation number to higher dimensions, especially
dimension two and homeomorphisms of the torus.

In this article we consider the rotation set $\cR(F)$ as defined in
[MZ], for an area preserving lift $F$ of a homeomorphism $f: \t^2\to
\t^2$, which is homotopic to the identity.  There are two fundamental
questions about this set, both of which are only partially resolved.
The first is what subsets of the plane can be realized as a rotation
set.  It is known that the set must be compact and convex and it is
known (see [K]) that any convex polygon with rational vertices can be
realized.  The second, more interesting, question is what is the
relation between this set and the existence of periodic points for
$f$.  For the case when $\cR(F)$ has interior in the plane, a great
deal is known about this question even without the hypothesis that $f$
is area preserving (see [ML] and [F2]).

The first of these questions seems quite difficult in the case that
$\cR(F)$ is an interval.  It is not known, for example if $\cR(F)$
can be a line segment containing no rational points (i.e. points in
$\r^2$ with both co-ordinates rational.)  In this article we focus on
the second question in the case when $\cR(F)$ has no interior in the
plane, i.e. is an interval.  We show that in this case if a vector $v$
lies in $\cR(F)$ and has both co-ordinates rational, then there is a
periodic point $x\in \t^2$ with the property that
$$\frac{F^q(x_0)-x_0}q = v$$ where $x_0\in \re^2$ is any lift of $x$ and
$q$ is the least period of $x$.

\ms

\heading \S1 Background and Definitions \endheading

Suppose $f: \t^2\to \t^2$ is a homeomorphism homotopic to the identity map, and 
let $F: \re^2\to \re^2$ be a lift.

\proclaim{(1.1) Definition} 
Let $\cR(F)$ denote the set of accumulation 
points of the subset of $\re^2$ $$\left\{ \frac{F^n(x)-x}n \biggm|
x\in \re^2, n\in Z^+ \right\},$$ thus $\nu\in\cR(F)$ if there are
sequences $x_i\in \re^2$ and $n_i\in Z^+$ with $\lim n_i=\infty$ such
that $$\lim_{i\to\infty} \frac{F^{n_i}(x_i)-x_i}{n_i} = \nu.$$ 
\endproclaim

In [MZ] the rotation set is defined for a map homotopic to the
identity (rather than a homeomorphism) $f: \t^n\to \t^n$. However, we
shall be concerned only with homeomorphisms of $\t^2$. In [MZ] it is
shown that for homeomorphisms of $\t^2$, $\cR(F)$ is convex.

An interesting and important question is to determine precisely which 
compact convex subsets of $\r^2$ can be the rotation set for a homeomorphism
of the torus.  A result of Kwapicz [K] shows that any convex polygon
with rational vertices can be realized  as a rotation set.

Closely related to this is the definition of the rotation vector of a 
point $x \in \r^2$ for the lift $F$.

\proclaim{(1.2) Definition}
Suppose $f: \t^2 \to \t^2$ is a homeomorphism 
which is isotopic to the identity map and let
$F:\r^2 \to \r^2$ be a lift of $f$.
The rotation vector of $x$ under $F$ is defined by 
$$\rho(x,F) = \lim_{n\to\infty} \frac{F^n(x)-x}{n},$$
if this limit exists.
\endproclaim \ms

Of course, $\rho(x,F)$ may not exist for a given $x$ (more about this
below), but it always exists for some choices of $x$.
In fact it is shown in [MZ] that the rotation set $\cR(F)$ is the
convex hull of the set of rotation vectors of points $x \in \r^2.$

In the case that $f$ preserves a measure $\mu$ we can show the 
existence of the rotation vector $\rho(x,F)$ for a large set of
values of $x$.  We consider the function $\phi(x) = F(x) - x$. 
This function is defined on $\r^2$, but is periodic in the sense
that $\phi(x +(m,n)) = \phi(x)$ for every $(m,n) \in \z^2,$ and
hence can also be considered as a function defined on $\t^2.$

It is then easy to see that 
$$\align\frac{1}{N}\sum_{n=0}^{N-1}\phi(f^n(x)) &=
\frac{1}{N}\sum_{n=0}^{N-1}(F^{n+1}(y) - F^n(y))\\
&\\
 &= \frac{(F^n(y) - y)}{N}.\endalign$$
So
$$\rho(x,F) = \lim_{N \to \infty} \frac{1}{N}\sum_{n=0}^{N-1}\phi(f^n(x)).$$

It follows from the Birkhoff ergodic theorem that this last limit
exists for a set of values of $x$ which has full measure with respect
to $\mu.$ Also from this formula (or from the fact that $F(x +(m,n)) =
F(x) + (m,n)$ ) it follows that $\rho(x + (m,n), F) = \rho(x, F)$ for
all $(m,n) \in \z^2$, and hence $\rho(x, F)$ can be consider to be a
function defined on $\t^2.$   Then another consequence of the Birkhoff ergodic
theorem is that the function $\rho(x,F)$ is integrable and 
$$\int \rho(x,F) d\mu = \int \phi(x) d\mu.$$ 

\bigskip
\proclaim{(1.3) Definition}
Suppose $f: \t^2 \to \t^2$ is a homeomorphism of the
surface $\t^2$ which is isotopic to the identity map and preserves
a measure $\mu$  and let $F:\r^2 \to \r^2$ be a lift of $f$.
The mean rotation vector of $F$ is
an element of $\r^2$ denoted $\rho_{\mu}(F)$, and is defined by 
$$\rho_{\mu}(F) = \int \rho(x,F) d\mu.$$
\endproclaim \ms

A key property of the mean rotation vector is the following well
known and easy result a proof of which can be found in [F3].

\proclaim{(1.4) Proposition} Suppose $f$ and $g$ are homeomorphisms
of $\t^2$ which are isotopic to the identity and have lifts $F$ and
$G$.  If $f$ and $g$  preserve a probability measure $\mu$.  Then
$$\rho_{\mu}(F \circ G) = \rho_{\mu}(F) + \rho_{\mu}(G).$$
\endproclaim \ms

In this article we focus on a key property of the rotation set -- namely,
its relation to the existence of periodic points.  A considerable amount
is known in the case that the rotation set has interior in the plane.
In [F2] the following result is shown.

\proclaim{(1.5) Theorem} Suppose $f: \t^2\to \t^2$ is a 
homeomorphism homotopic to the identity and $F: \re^2\to \re^2$ is a
lift. If $\nu$ is a vector with rational co-ordinates in the interior
of $\cR(F)$, then there is a point $p\in \re^2$ such that $\pi(p)\in
\t^2$ is a periodic point for $f$ and 
$$\nu = \lim_{n\to\infty}\frac{F^n(p)-p}n.$$
\endproclaim

This result should be compared with the classical result of Poincar\'e
for homeomorphisms of the circle which says that the rotation number
is the rational in lowest terms, $p/q$, if and only if there exists a
periodic point of period $q$ with rotation number $p/q.$ 

The theorem about torus homeomorphisms naturally raises two questions:
what happens for rational points on the boundary and what happens when
$\cR(F)$ has no interior in $\r^2$?  The first of these questions is
resolved by an unpublished example due to Matt Grayson.  It shows that
there are homeomorphisms (which can even be chosen to be area
preserving) with $\cR(F)$ a polygon one edge of which contains
infinitely many rational points and no rational point in this side,
except the vertices, corresponds to a periodic point. It is shown in
[MZ] that a rational extreme point (in fact any extreme point) must
correspond to a periodic point.

Our main result in this article is to address the more difficult
question of what happens when the rotation set $\cR(F)$ has no
interior in the plane.  In this case it must be a line segment or
a single point.  We are able to resolve this question in the case
that the homeomorphism is area preserving.

The case that $f$ is area preserving and $\cR(F)$ consists of a 
single rational point is handled as a special case of a result
of [F1].  This result says that if the mean rotation vector of
$F$ is rational then there is a periodic point whose rotation vector
is equal to that vector.  If $\cR(F)$ consists of a single rational
vector then it is easy to see that the mean rotation vector must be
this vector.

Our main result, which deals with the case $\cR(F)$ is an interval,
is the following.

\proclaim{Theorem} Suppose $f: \t^2\to \t^2$ is an area preserving
homeomorphism homotopic to the identity and $F: \re^2\to \re^2$ is a
lift. If $\cR(F)$ is an interval and $\nu$ is a vector with rational
co-ordinates in it, then there is a point $p\in \re^2$ such that
$\pi(p)\in
\t^2$ is a periodic point for $f$ and
$$\nu = \lim_{n\to\infty}\frac{F^n(p)-p}n.$$
\endproclaim

There are several interesting unanswered questions related to this
theorem.  First, it is not known whether the hypothesis that $f$ be
area preserving is necessary.  Second it is not known which intervals
can occur as the rotation set of a homeomorphism of $\t^2.$  It is 
easy to see that any interval with rational slope can be realized.
An example due to Katok and discussed in [FM] shows that
an interval with one endpoint rational and no other rational points
can be realized.  But it is not known whether an interval with 
no rational points or a single interior rational point can be 
realized.
\bigskip

A key ingredient in the proof of our main theorem is the concept of
$\varep$-chain.

\proclaim{(1.6) Definition} An {\it $\varep$-chain\/} for $f$ is a sequence 
$x_1,x_2,\dots,x_n$ of points in $X$ such that
$$d(f(x_i), x_{i+1})<\varep \qquad\text{for } 1\le i\le n-1.$$
If $x_1=x_n$ it is called a periodic $\varep$-chain. 
\endproclaim

The relationship of $\varep$-chains to our program of finding 
periodic points is given by Lemma (1.7) below.  A proof of this
lemma can be found as (2.1) of [F2].
We assume  that $f: \t^2\to \t^2$ is a homeomorphism homotopic
to the identity and $F: \re^2\to \re^2$ is a lift, i.e., if $\pi:
\re^2\to \t^2$ is the covering projection then $\pi\circ F=f\circ\pi$.

\proclaim{(1.7) Lemma} If $F$ has no fixed points, then there is an $\varep>0$ 
such that no periodic $\varep$-chain for $F$ exists. \endproclaim

\bigskip

\heading \S2 Proof of the Main theorem \endheading

In this section we consider area preserving homeomorphisms of the 
torus, by which we mean a homeomorphism $f$ which preserves Lebesgue
measure $\mu$.   This will be part of the hypothesis of our main
result Theorem (2.2) below.  However, it should be observed that
it is an immediate consequence of this result that the same theorem
is valid for any measure $\mu_0$ which is homeomorphic to Lebesgue
measure in the sense that there is a homeomorphism $h: \t^2 \to \t^2$
with the property that $\mu_0(h(B)) = \mu(B)$ for every Borel set
$B$ in $\t^2.$  This is because a homeomorphism $f_0$ will preserve
$\mu_0$ if and only if the conjugate homeomorphism $f = h^{-1} \circ f_0
\circ h$ preserves Lebesgue measure $\mu.$

Theorem (2.2) actually a corollary of the following

\proclaim{(2.1) Proposition}  Suppose $f: \t^2\to \t^2$ is an area
preserving homeomorphism homotopic to the identity and $F: \re^2\to
\re^2$ is a lift. Let $K_\varep$ denote the convex hull of $\cR(F)$
and the ball of radius $\varep$ centered at $\rho_\mu(F)$.  If $\nu$
is a vector with rational co-ordinates in the interior of $K_\varep$,
for all $\varep > 0$, then there is a point $z\in \re^2$ such that
$\pi(z)\in \t^2$ is a periodic point for $f$ and 
$\nu = \rho(z,F)$.  Moreover, if $\nu = (p/q, r/q)$ where $p,\ q$, and
$r$ are relatively prime then $q$ is the period of $\pi(z)$.
\endproclaim

Before proving this we show that our main result follows from this.

\proclaim{(2.2) Theorem} Suppose $f: \t^2\to \t^2$ is an area preserving
homeomorphism homotopic to the identity and $F: \re^2\to \re^2$ is a
lift. If $\cR(F)$ is an interval and $\nu$ is a vector in it with rational
co-ordinates, then there is a point $z\in \re^2$ such that
$\pi(z)\in
\t^2$ is a periodic point for $f$ and
$\nu = \rho(z,F)$.  Moreover, if $\nu = (p/q, r/q)$ where $p,\ q$, and
$r$ are relatively prime then $q$ is the period of $\pi(z)$.
\endproclaim

\demo{Proof} Since $\rho_\mu(F) = \int \rho(x, F) d\mu$ and 
$\cR(F)$ is the convex hull of the set of vectors $\{\rho(x, F)\}$,
it must be the case that  $\rho_\mu(F)$ is in the interval
$\cR(F)$.  Therefore if $\nu$ is in $\cR(F)$ and not an endpoint,
it follows that $\nu \in K_\varep$ for all $\varep > 0$.  Hence
in this case the desired $z$ exists by (2.1).  

If $\nu$ is an endpoint of $\cR(F)$ then it is an extreme point and
results of [MZ] show that $\nu$ is the rotation vector of an ergodic
invariant measure for $f$.  In this case it is an easy consequence of
a result of [F1] (Theorem (3.5 )) that the desired $z$ exists.
\hfill\bbox 
\enddemo

We turn now to the proof of proposition above, beginning with
two lemmas.

\proclaim{(2.3) Lemma} Suppose $f: \t^2\to \t^2$ is as in (2.1) and
$v \ne 0$ is one of the end points of the interval $\cR(F)$.  Then for any
$\varep > 0$ and any $\delta >0$ there exists an $\varep$-chain for $F$ from
$0 \in \r^2$ to $w \in \z^2$ with the property that
$$\left\| \frac{v}{ \|v\|} - \frac{w}{\|w\|} \right\| < \delta.$$
\endproclaim
\bigskip

\demo{Proof}
Because $v$ is an extreme point of $\cR(F)$ it follows that there
is an ergodic invariant measure with mean rotation vector $v$ (see
[MZ]).  From this is follows that there is a point $x_0 \in \r^2$ with
$\rho(x_0, F) = v$.  Of course this implies that $\rho(x_0 + m, F) = v$,
for any $m \in \z^2$.  Moreover, by the Poincar\'e recurrence theorem
we can assume that $\pi(x_0)$ is a recurrent point of $f:\t^2 \to \t^2.$
Thus there infinitely many $n$ such that $\pi(F^n(x_0))$ is within
$\varep$ of $\pi(x_0)$.  Hence we can choose arbitrarily large $n$
such that the sequence $x_0, F(x_0), F^2(x_0), \dots F^{n-1}(x_0),
x_0 + v(n)$ is an $\varep$-chain from $x_0$ to $x_0 + v(n)$ for
some $v(n) \in \z^2$.  In addition $v(n)$ has two important properties:
$$\lim_{n \to \infty} \|v(n)\| = \infty, \tag{1}$$
and
$$\lim_{n \to \infty} \frac{v(n)}{\|v(n)\|} = \frac{v}{\|v\|}. \tag{2}$$

Since $f$ preserves area on $\t^2$ it follows
that a set of full (Lebesgue) measure in $\t^2$ consists of recurrent points.
From this it is easy to see that $f$ is chain transitive, (see [F1]) and
in particular that there is an $\varep$-chain from $0 \in \t^2$ to 
$\pi(x_0)$ and also one from $\pi(x_0)$ to $0.$  Lifting the first of these to
an $\varep$-chain for $F$ starting at $0$ we obtain a  $\varep$-chain from
$0$ to $x_0 + w_0$ for some $w_0 \in \z^2.$  A translate by $w_0$ of the
$\varep$-chain constructed above will be an $\varep$-chain from $x_0 + w_0$
to $x_0 + v(n) + w_0$.  Concatenating gives an $\varep$-chain from $0$
to $x_0 + v(n) + w_0$.  Starting with the $\varep$-chain from $\pi(x_0)$
to $0$ and lifting, translating and concatenating we get an $\varep$-chain
from $0$ to $x_0 + w_0$ to $x_0 + w_0 +v(n)$ to $w_0 + v(n) +w_1$ for
some $w_1 \in \z^2$.  The integer vectors $w_0$ and $w_1$ are independent
of $n$.  Thus it follows from (1) and (2) above that
$$\lim_{n \to \infty} \frac{w_0 + v(n) + w_1}{\|v(n)\|} = \frac{v}{\|v\|}.$$
Hence for $n$ sufficiently large, the $\varep$-chain from $0$ to
$w = w_0 +v(n) +w_1$ will serve as the desired one.
\hfill\bbox 
\enddemo
\bigskip

\proclaim{(2.4) Lemma} Suppose $f: \t^2\to \t^2$ is as in (2.1) and
$v \ne 0$ is in the ball of radius $\varep/2$ centered at $\rho_\mu(F)$.
Then for any $\delta >0$ there exists an
$\varep$-chain for $F$ from $0 \in \r^2$ to $w \in \z^2$ with the
property that $$\left\| \frac{v}{ \|v\|} - \frac{w}{\|w\|} \right\| <
\delta.$$
\endproclaim

\demo{Proof}
Let $u = v - \rho_\mu(F)$ and consider the homeomorphism $G: \r^2 \to
\r^2$ given by $G_0(x) = F(x) + u.$  Clearly $G_0$ is a lift of a 
homeomorphism $g_0$ of $\t^2$ and from (1.4) it follows that 
$$\rho_\mu(G_0) = \rho_\mu(F) + u = v.$$

The Oxtoby--Ulam Theorem (see [OU]) asserts that the homeomorphism $g_0$ can
be perturbed by an arbitrarily small amount in the $C^0$ topology
to obtain a homeomorphism $g$ which is ergodic with respect to
Lebesgue measure.  By making $g$ sufficiently close to $g_0$ and choosing
the lift $G$ of $g$ which is close to $G_0$ we can assume that
$\|\rho_\mu(G) - \rho_\mu(G_0)\|$ is as small as we wish.  In particular
since $v =\rho_\mu(G_0)$ we can assume that
$$\left\|\frac{\rho_\mu(G)}{\|\rho_\mu(G)\|} - 
\frac{v}{\|v\|}\right\| < \delta.$$

The fact that $g$ is ergodic implies that for all $x$ in $\r^2$ except
for a set of measure zero we have that $\pi(x)$ is recurrent under $g$,
and $\rho(x, G) = \rho_\mu(G).$  Conjugating $G$ by a very small translation
(and continuing to call the new homeomorphism $G$) we can assume that
$0\in\r^2$ is in this set of full measure.  Hence we can assume that
$\pi(0)$ is recurrent under the homeomorphism $g$ and
$$\left\|\frac{\rho(0,G)}{\|\rho(0,G)\|} - 
\frac{v}{\|v\|}\right\| < \delta/2.$$

It follows that there are large values of $n$ such that
$G^n(0)/n$ is as close as we wish to $\rho(0,G)$ and $G^n(0)$ is as close
as we wish to an element $w(n) \in \z^2.$  For such an $n$ we note
that the orbit segment (of $G$),  $0, G(0), G^2(0), \dots, G^n(0)$ is
an $\varep$-chain for $F$ and if we replace $G^n(0)$ by $w(n)$ it is
a $\varep$-chain from $0$ to $w = w(n)$.  Also if $n$ has been chosen
sufficiently large
$$\left\|\frac{w}{\|w\|} - 
\frac{v}{\|v\|}\right\| < \delta.$$
\hfill\bbox 
\enddemo
\bigskip

We return now to the proof of the main result.

\demo{Proof of (2.1)} We note first that it suffices to consider a special
case; namely it is enough to consider the case when $\nu = 0$ and prove
that in this case there is a fixed point for the homeomorphism
$F: \r^2 \to \r^2$.  This is  because if $\nu = (p/q, r/q)$ where $p,\ q$, and
$r$ are relatively prime then we can consider the new homeomorphism
$G = F^q - (p,r)$.  $G$ is a lift of $f^q$ and a fixed point $z$ of $G$
in $\r^2$ will satisfy $\rho(z,G) = 0$ from which it follows that
$0 = \rho(z,G) = \rho(z, F^q - (p,r)) = \rho(z, F^q) - (p,r) = 
q \rho(z,F) - (p,r)$. So $\rho(z,F) = (p/q, r/q) = \nu$.  Moreover,
$z$ is a fixed point of $G$ so $\pi(z)$ is a fixed point of $f^q$.
If the period of $z$ were less than $q$ then  $\rho(z,F)$ could
be written with a denominator less than $q$ which would contradict
the assumption that  $p,\ q$, and $r$ are relatively prime.

Thus we may assume that $\nu = 0$ and we must prove that $F$ has
a fixed point $z$.  Our strategy is to construct a periodic $\varep$-chain
for $F$ for each $\varep > 0$.  It then will follow from (1.7) that 
$F$ has a fixed point. 

Since $\nu$ is in the interior of the convex set $K_\varep$, there are 
finitely many points $v_i$ such that $\nu$ is in the convex hull 
of $\{v_i\}$ (see Steinitz's theorem in [HDK]).  Clearly we can assume each 
$v_i$ is either an endpoint of $\cR(F)$ or a point in the ball of
radius $\varep$ centered at $\rho_\mu(F).$

Thus by (2.3) and (2.4) above, for any $\delta >0$ we can find for 
each $v_i$ a corresponding $w_i \in \z^2$ with the property that there
is a $\varep$-chain for $F$ from $0$ to $w_i$ and 
$$\left\|\frac{w_i}{\|w_i\|} - 
\frac{v_i}{\|v_i\|}\right\| < \delta.$$

The fact that $0$ is in the interior of the convex hull of $\{v_i\}$
implies that if $\delta$ is chosen sufficiently small we can find
positive rational numbers $A_i$ such that 
$$\sum_i A_i w_i = 0.$$
Clearing fractions we may assume that each $A_i$ is a positive integer.

If we now concatenate the $\varep$-chain from $0$ to $w_1$ with
a translate of itself we get an $\varep$-chain from $0$ to $w_1$
to $2w_1$.  Repeating this until we have concatenated  $A_1$ 
of these $\varep$-chains we get an $\varep$-chain from $0$ to
$A_1 w_1$.  Starting at this point we concatenate $A_2$ translates
of the $\varep$-chain from $0$ to $w_2$ and then repeat this for
$w_3$ etc.  We finish with an $\varep$-chain from $0$ to 
$\sum A_i w_i = 0,$ i.e. a periodic $\varep$-chain.  As we remarked
above the fact that this can be done for any $\varep > 0$ implies
by (1.7) that $F$ has a fixed point.  This completes the proof of
(2.1).

\hfill\bbox 
\enddemo

\enddocument

\bye